\documentclass[11pt]{article}

\usepackage{amsmath}
\usepackage{amssymb}
\usepackage{latexsym}

\def\cF{{\cal F}} 
\def\vf{\varphi}
\def\C{{\cal C}}
\def\D{\Delta}
\def\wt{\widetilde}
\def\cD{{\cal D}}
\def\Supp{\rm Supp }
\def\eqdefa{\buildrel\hbox{\footnotesize def}\over =}
\def\longformule#1#2{
\displaylines{
\qquad{#1}
\hfill\cr
\hfill {#2}
\qquad\cr
}
}
\def\norm#1#2{\|#1\|_{#2}}
\def\normsup#1{\|#1\|_{L^\infty}}

\def\refeq#1{~(\ref{#1})}
\def\ccite#1{~\cite{#1}}

\newcommand{\RR}{{\Bbb R}}

\newcommand{\ZZ}{{\Bbb Z}}
\newcommand{\NN}{{\Bbb N}}

\newcommand{\Sum}{\displaystyle \sum}
\newcommand{\Sup}{\displaystyle \sup}

\def\eps{\varepsilon}
\def\s{\sigma}
\def\b{\beta}
\def\a{\alpha}
\def\lam{\lambda}
\renewcommand{\div}{{\rm div}}

\def\p{\partial}

\newtheorem{theo}{Theorem}[section]

\newtheorem{lem}[theo]{Lemma}

\newtheorem{defi}[theo]{Definition}
\newtheorem{cor}[theo]{Corollary}

\begin{document}

\title{ Global well-posedness for a Smoluchowski equation coupled with Navier-Stokes equations in $2$D. }
\author
{ P. Constantin\\
Department of Mathematics, The University of Chicago\\ 5734 S. University Avenue, Chicago, Il 60637\\email const@math.uchicago.edu\\and \\
\\
Nader Masmoudi \\
 Courant Institute, New York University \\
  251 Mercer St,
New York NY 10012\\
email:masmoudi@cims.nyu.edu }

%\date{}

\maketitle

\vspace{8mm}

\noindent \underline{\bf Abstract}

{\sl We prove global existence for a  nonlinear Smoluchowski equation (a nonlinear Fokker-Planck equation) coupled with Navier-Stokes equations in $2$d. The proof  uses a deteriorating regularity estimate in the spirit of \cite{CM01} (see also \cite{BC94})}



\section{Introduction}

Systems coupling fluids and particles are of great interest in many branches of applied physics and chemistry. The equations attempt to describe the behavior of complex mixtures of particles and fluids, and as such, they present numerous challenges, simultaneously at three levels: at the level of 
their derivation, the level of their numerical simulation and that of their mathematical  treatment. In this paper we concentrate solely on one 
aspect of the mathematical treatment, the regularity of solutions. The particles in the system are described by a probability distribution
$f(t,x,m)$ that depends on time $t$, macroscopic variable $x\in {\mathbb R}^n$,
and particle configuration $m\in M$. Here $M$ is a smooth compact Riemannian
manifold without boundary. The particles are transported by a fluid, agitated by thermal noise, and interact among themselves. This is reflected in a kinetic
equation for the evolution of the probability distribution of the particles 
(\cite{BCAH87,doi}). The interaction between particles -- a micro-micro interaction -- is modeled in a mean-field fashion by a potential that  represents the tendency of particles to favor certain coherent configurations. The interaction between particles occurs only when the concentration of particles is sufficiently high. Mathematically, this term is responsible for the nonlinearity of the Smoluchowski (Fokker-Planck) equation, and physically, it is responsible for nematic phase transitions. Because the particles are considerably small, and for smooth flows, the
Lagrangian transport of the particles is modeled using a Taylor expansion of the velocity field. This gives rise to a drift term in the Smoluchowski equation that depends on the spatial gradient of velocity. It is a macro-micro term, and it causes mathematical difficulties in the regularity theory.

The fluid is described by the incompressible Navier-Stokes equations. The
microscopic particles add stresses to the  fluid. This is the micro-macro
interaction and it is the most puzzling and important physical aspect of the 
problem. Indeed,  while a
macro-micro interaction can be derived, in principle,  by assuming that the macroscopic entities vary little on the scale of the microscopic ones, the ``scaling up'' of the effect of microscopic quantities to the macroscopic level is more mysterious. A principle based on an energy dissipation balance, and 
that recovers familiar results in simple cases was proposed in \cite{const}, where the regularity of  nonlinear Fokker-Planck systems coupled with Stokes equations in 3D was also proved. The linear
Fokker-Planck system coupled with Stokes equations was considered in \cite{otto}. The nonlinear Fokker-Planck equation driven by a time averaged 
Navier-Stokes system in 2D was studied in \cite{CFTZ06}.

An approximate closure of the linear Fokker-Planck equation reduces the 
description to closed viscoelastic equations for the 
added stresses themselves. This leads to well-known  non-Newtonian fluid models that have been studied extensively.  
For regularity results we refer to Lions and Masmoudi  \cite{LM00cam} where the 
existence of global weak solutions was proved for an 
Oldroyd-type model. 
In  Guillop{\'e}  and Saut \cite{GS90} and  \cite{GS90b}, 
the existence of local strong solution was proved. Also, 
Fern{\'a}ndez-Cara,  Guill{\'e}n  and
             Ortega  \cite{FGO98},   \cite{FGO97} and \cite{FGOb}
proved local well posedness in Sobolev spaces. 
We also mention Lin, Liu and  Zhang  \cite{LLZ05} where 
a formulation based on the deformation tensor is used to 
study the Oldroyd-B model. 

An other model for the polymers is the FENE dumbbell model. 
From mathematical point of view, this model was studied by several authors.
In particular W. E, Li and Zhang \cite{ELZ04},  
Jourdain, Lelievre and Le Bris \cite{JLL04} and  
   Zhang and Zhang \cite{ZZ06arma}
 proved local well-posedness.  Moreover,  Lin, Liu and  Zhang  \cite{LLZ07}
  proved global existence near equilibrium.

% We also mention that there is an other 
% micro-macro 
% models were considered. 

\subsection{The model}
Consider the system
\begin{equation}
\label{NS}
 \ \left\{
\begin{array}{rclll}
\frac{\partial v}{\partial t} +   v\cdot\nabla  v - \nu
    \Delta v +
\nabla p & =
&    \nabla\cdot \tau   & \mbox{ in }&\Omega\times (0,T)\\
   \\
\frac{\partial f}{\partial t} +   v\cdot\nabla  f+   \div_g (G(v,f)  f  ) - \Delta_g f    & = & 0  &\mbox{ in
}&\Omega\times (0,T)\\
\\ 
\div  v &=&0&\mbox{ in }&  \Omega\times(0,T), 
\end{array}
\right.
\end{equation}
where $ \tau_{ij} = \int_M  \gamma^{(1)}_{ij}(m)  f(t,x,m) dm   +
 \int_M  \int_M  \gamma^{(2)}_{ij}(m_1,m_2)  f(t,x,m_1) f(t,x,m_2)   dm    $. 
We denote  $G(v,f) = \nabla_g U + W $ 
where  $W = c^{ij}_\a \partial_j v_i $ and   $U = K f$ is a potential 
  given by  
\begin{equation}
U(t,x,m) = \int_M K(m,q) f(t,x,q)\, dq
\end{equation}
with a kernel $K$ which is a smooth, time and space  independent symmetric 
function $K : M\times M \to \RR$.  We also take $\Omega = \RR^2$. 

\subsection{Statement of the result}

\begin{theo} \label{main} 
Take  $v(0) \in W^{1+\eps_0,r}\cap L^2 (\RR^2)$ and 
$f(0) \in W^{1,r} (H^{-s})$, for some $r>2$ and 
$\eps_0 > 0$  and  $f \geq 0$, 
   $      \int_M f_0  \in L^1 \cap L^\infty$. Then (\ref{NS}) has 
a global solution in 
$ v \in L^\infty_{loc} ( W^{1,r} ) \cap L^2_{loc} ( W^{2,r}) $ and 
$ f \in  L^\infty_{loc} ( W^{1,r} (H^{-s})  ) $. 
Moreover, for $T>T_0 > 0$, we have $v \in  L^\infty( (T_0,T); W^{2-\eps,r} )$.  
\end{theo}

\subsection{Preliminaries}
We  define
$\C$ to be the ring of center
$0$, of small radius $1/2$ and great radius $2$. There exist two
nonnegative  radial
functions $\chi$ and $\vf$ belonging respectively to~${\cal D} 
(B(0,1)) $ and to
${\cal D} (\C) $ so that
\begin{equation}
\label{lpfond1}
\chi(\xi) + \sum_{q\geq 0} \vf (2^{-q}\xi) = 1,
\end{equation}
\begin{equation}
\label{lpfond2}
|p-q|\geq 2
\Rightarrow
\Supp\ \vf(2^{-q}\cdot)\cap \Supp\ \vf(2^{-p}\cdot)=\emptyset.
\end{equation}

\noindent
For instance, one can take $\chi \in \cD (B(0,1))$ such that $
\chi  \equiv 1 $ on $B(0,1/2)$ and take
$$
\vf(\xi) = \chi(\xi/2) -
\chi(\xi).
$$
Then, we are able to define the Littlewood-Paley decomposition. Let us denote
by~$\cF$ the Fourier transform on~$\RR^d$. Let
$h,\
\wt h,\  \D_q,\  S_q$ ($q \in \ZZ$) be defined as follows:
%\noindent{\bf Notations}

$$\displaylines{
\label{defnotationdyadique}h = {\cal F}^{-1}\vf\quad {\rm and}\quad \wt h =
{\cal F}^{-1}\chi, \cr
\D_qu = \cF^{-1}(\vf(2^{-q}\xi)\cF u) = 2^{qd}\int h(2^qy)u(x-y)dy,\cr
S_qu
%= \sum_{p\leq q-1} \D_pu
=\cF^{-1}(\chi(2^{-q}\xi)\cF u) =2^{qd} \int \wt h(2^qy)u(x-y)dy.\cr
}
$$

We  use  the  para-product decomposition of Bony 
(\cite{Bony81})
% \cite{Bony82}
$$
uv = T_uv \ + \ T_vu\ + \ R(u,v)
$$
where
$$
T_uv  =  \Sum_{q\in \ZZ} S_{q-1}u \Delta_q v \quad\hbox{and}
\quad   R(u,v)   = \Sum_{|q-q'|\leq 1 } \Delta_{q'}u \D_qv.
$$

We define the inhomogeneous and  homogeneous Besov spaces by 

\begin{defi}
Let $s$ be a real number, p and r two real numbers greater than~$1$.
Then we define the following norm
$$
\|u\|_{\wt B^s_{p,r}} \eqdefa \|S_0u\|_{L^p}+\Bigl\|\left( 2^{qs} \|\Delta
_qu\|_{L^p}\right)_{q\in \NN}\Bigr\|_{\ell^r(\NN)}
$$
and the following semi-norm
$$
\|u\|_{B^s_{p,r}} \eqdefa\Bigl\|\left( 2^{qs} \|\Delta
_qu\|_{L^p}\right)_{q\in \ZZ}\Bigr\|_{\ell^r(\ZZ)}.
$$
\end{defi}

\begin{defi}${\atop}$
\begin{itemize}
\item
Let $s$ be a real number, p and r two real numbers greater than~$1$.
We denote by~$\wt B^s_{p,r}$ the space of tempered distributions~$u$ such
that~$\|u\|_{\wt B^s_{p,r}}$ is finite.
\item
If $s < d/p$ or $s = d/p$ and $r=1$ we define
the homogeneous Besov space~$B^s_{p,r}$ as  the closure of compactly supported
   smooth functions for  the norm~$\|\cdot\|_{B^s_{p,r}}$.
\end{itemize}
\end{defi}

We refer to~\cite{Chemin99}  for the
proof of the following results and for the multiplication law
in Besov spaces.
\begin{lem}
\label{loc}
    $$  \| \D_q u \|_{L^b} \leq 2^{d({1 \over a} -{1 \over b})q } \|\D_q u
\|_{L^a}
\quad{\rm for }\ b\geq a \geq 1
$$
    $$ \|  e^{t \D} \D_q u \|_{L^b} \leq C 2^{-ct 2^{2q}} \| \D_q u
\|_{L^b}
$$
\end{lem}
The following corollary is straightforward.
\begin{cor}
\label{inclusiontypesob}
{\sl
If~$b\geq a\geq 1$, then, we have the following continuous embeddings
\[  
B^s_{a,r} \subset
B^{s-d\Bigl( \frac 1 a-\frac 1 b\Bigr)}_{b,r}.
\]
}
\end{cor}

\begin{defi}
Let~$p$ be in~$[1,\infty]$ and~$r$ in~$\RR$; the space $ \widetilde
L^p_T(C^r) $ is
the space of  distributions~$u$ such that
$$
\|u\|_{\widetilde L^p(0,T;C^r)}\eqdefa \Sup_q  2^{qr} \| \D_q
u\|_{L^p_T(L^\infty)} <
\infty.
$$
\end{defi}
We will use the following theorem from \cite{CM01}

\begin{theo}
\label{estimation2dprecisee}
Let $v$ be the solution in~$L^2_T(H^1)$
  of the two dimensional Navier-Stokes system
$$
(NS_\nu)
\left\{
\begin{array}{rcl}
\displaystyle  \frac {\partial v}{\partial t}+v\cdot\nabla v -\nu\Delta v & =
&-\nabla p+ f
\\
\div  v & = & 0 \\
v|_{t=0} & = & v_0.
\end{array}
\right.
$$
with an initial data
in~$L^2$  and an external force
$f$ in~$L^1_T(C^{-1})\cap L^2_T(H^{-1})$; then, for any  $\eps$,
a~$T_0$ in the interval~$]0,T[$ exists such  that
\[
\|\nabla v\|_{\wt L^1_{[T_0,T]}(C^0)}\leq \eps.
\]
\end{theo}

\section {A deteriorating regularity estimate}
\label{estimationapertetypebc}
The main part of this section is the proof of a deteriorating regularity 
estimate for transport equations in the 
spirit of \cite{BC94} and \cite{CM01}.
 After this proof, we will apply this
estimate in order to prove Theorem \ref{main}.

We also denote $ H = (-\Delta_g + I  )^{-s/2}$  with $s> d/2+1$. 

\begin{theo}
\label{besovinegapertetransport}
Let~$\sigma$ and~$\beta$ be two elements 
 of~$]0,1[$ such that $\s+\b<1$. A constant~$C$ exists that
satisfies the
following properties.
Let~$T$ and~$\lambda$ be two positive numbers and~$v$ a smooth divergence
free vector field so that
\begin{equation} 
\label{besovinegapertetransporteq1}
\sigma-\lambda\|\nabla v\|_{\wt L^1_T(C^0)}\geq \b   
\quad {\rm and} \quad  \sigma+\lambda\|\nabla v\|_{\wt L^1_T(C^0)} \leq 1-\b .
\end{equation}
Consider two smooth functions~$f$ and~$v$ so that~$f$ is the solution of
\begin{equation} \label{transp}
\left\{\begin{array}{rcl}
\partial_tf+v\cdot\nabla f + \div_g (G(v,f)  f  ) - \Delta_g f    & = & 0\\
f_{|t=0} & = & f_0.
\end{array}
\right.
 \end{equation}

Then we have, if~$\lam\geq 3C$,
\begin{equation}
M^\sigma_\lambda(f)  \leq  3 \|f_0\|_{B^\sigma_{p,\infty} (H^{-s})   } 
 + \frac {3C}\lam 
M_\lam^{\s+1}(v) \\
\end{equation}
where
\begin{eqnarray} 
\label{definregulariteaperte}
M^\s_\lam(v) & \eqdefa &
\Sup_{t\in[0,T],q}  
2^{q\s-\Phi_{q,\lam}(t)}\|\D_qv(t)\|_{L^p}\quad\hbox{or}\\
M^\s_\lam(f) & \eqdefa &
\Sup_{t\in[0,T],q}  
2^{q\s-\Phi_{q,\lam}(t)}\|\D_qf(t)\|_{L^p(H^{-s}) }\quad\hbox{with}\\
\label{definphasegronwall}
   \Phi_{q,\lam}(t,t') & \eqdefa & \lam\int_{t'}^t (\|S_{q-1}\nabla
v(t'')\|_{L^\infty} + 1)  dt''  ,\ \Phi_{q,\lam}(t) =
\Phi_{q,\lam}(t,0).
\end{eqnarray}
\end{theo}

We will use the notation $f_q\eqdefa
\Delta_qf$. Applying the operator $\Delta_q$ to 
the transport equation (\ref{transp}), we get 
 
\begin{equation} \label{transp-q}
\left\{\begin{array}{rcl} 
 \partial_tf_q+S_{q-1}v\cdot\nabla f_q +  \div_g (G(S_{q-1} v,S_{q-1}f)  f_q  ) 
 - \Delta_g f_q    +  R_q(v,f) = 0\\
{f_q}_{|t=0}  =   \D_qf_0.
\end{array}
\right.
\end{equation}
where $R_q$ is a rest term. 

We denote 
\begin{equation}
N_q^2(t,x) =  \int_M |H f_q |^2  dm 
\end{equation}
Applying $H$ to (\ref{transp-q}) and taking the $L^2$ norm on $M$, 
we get

\begin{equation} \label{transp-N-q} 
 \partial_t N_q^2+S_{q-1}v\cdot\nabla N_q^2 +  V (S_{q-1}v, S_{q-1}f,  f_q)    
  +   |\nabla_g H  f_q|^2      +    \int_M H f_q (H R_q(v,f)) dm  = 0 
\end{equation}
where 
   
\begin{equation} \label{def-V}
 V (v, h,   f) =  \partial_j v_i   \int_M  (H \div_g (c^{ij}_\a  f  )) (Hf) dm. 
 +  \int_M  (H \div_g ( \nabla_g h  f   )) (Hf) dm. 
\end{equation}  
 Hence, arguing as in \cite{CFTZ06}, we have 
 $| V (S_{q-1} v,  f_q)   | \leq C  ( | \nabla S_{q-1} v| + || S_{q-1} f  ||_{L^2(M)} ) N_q^2$.

We will use now the following lemma, postponing its proof: 
\begin{lem}
$R_q(v,f)$ satisfies 
\begin{eqnarray}
2^{q\s-\Phi_{q,\lam}(t)}\|H R_q(v(t),f(t))\|_{L^p(L^2)} \leq
Ce^{C\lam\|\nabla v\|_{\wt L^1_T(C^0)}} & & \nonumber\\
\label{besovinegapertetransportlemmeeq7}
   & &
\!\!\!\!\!\!\!\!\!\!\!\!\!\!\!\!\!\!\!\!\!\!\!\!\!\!\!\!\!\!\!\!\!\!\!\!\!\!
\!\!\!\!\!\!\!\!\!\!\!\!\!\!\!\!\!\!\!\!\!\!\!\!\!\!\!\!\!\!\!\!\!
\!\!\!\!\!\!\!\!\!\!\!\!\!\!\!\!\!\!\!\!\!\!\!\!\!\!\!\!\!\!\!\!\!\!\!\!\!\!
\!\!\!\!\!\!
{}\times \biggl(    
% \|\nabla_g H f(t)\|_{L^\infty (L^2)}    
M^{\s+1}_\lam(v)
+\Bigl( 1+ \|S_q\nabla v(t)\|_{L^\infty} + \sum_{|q'-q|\leq N}\| \D_{q'}\nabla
v(t)\|_{L^\infty}\Bigr)M^\s_\lam(f)\biggr).
\end{eqnarray}
\end{lem}
  
% Notice that since $\int |f| dm \leq 1$, we deduce that 
% $ \| H f \|_{L^\infty (\RR^2, L^2( M))} \leq C $. 

Taking the $L^p$ norm of $N_q$, we get 

\[  
\|N_q(t)\|_{L^p} \leq \|N_q(0)\|_{L^p} 
+\int_0^t\|H R_q(v(t'),f(t'))\|_{L^p(L^2)}  +   \| \nabla S_q v (t')\|_{L^\infty} 
  \|N_q(t')\|_{L^p}   dt'.
\]

After multiplication by~$2^{q\s-\Phi_{q,\lam}(t)}$,  we get

\begin{eqnarray*}
2^{q\s-\Phi_{q,\lam}(t)}\|N_q(t)\|_{L^p} \leq
2^{q\s}\|N_q(0)\|_{L^p} &+&
   \int_0^t
2^{-\Phi_{q,\lam}(t,t')}2^{q\s-\Phi_{q,\lam}(t')}   \| \nabla S_q v (t')\|_{L^\infty} 
    \| N_q \|_{L^p}dt' \\
&+& \int_0^t 2^{-\Phi_{q,\lam}(t,t')}
2^{q\s-\Phi_{q,\lam}(t')}\|H R_q(v(t'),f(t'))\|_{L^p(L^2)}dt'.
\end{eqnarray*}

Then, using the inequality (\ref{besovinegapertetransportlemmeeq7}) and taking 
the sup over $q$, we get

\begin{eqnarray}
M_\lam^\s(f) \leq
\|f_0\|_{B^\s_{p,\infty} (H^{-s})} + e^{C\lam\|\nabla v\|_{\wt
L^1_T(C^0)}} 
\Sup_{t\in[0,T],q}    \int_0^t 2^{-\Phi_{q,\lam}(t,t')} \\
\times \biggl(M_\lam^{\s+1}(v) + M_\lam^\s(f) \Bigl(1+ 
2 \|S_q\nabla v(t')\|_{L^\infty}
+\sum_{|q'-q|\leq N}\| \D_{q'}\nabla v(t')\|_{L^\infty}\Bigr)\biggr)dt'.
\end{eqnarray}

As~$\lam\|\nabla v\|_{\wt L^1_T(C^0)}$ is smaller than~$(\s-\b)$, we have
\[
e^{C\lam\|\nabla v\|_{\wt L^1_T(C^0)}}\leq e^{C(\s-\b)}.
\]
Moreover, by
definition of~$\Phi_{q,\lam}(t,t')$, it is obvious that
$$
% \int_0^t 2^{-\Phi_{q,\lam}(t,t')}\|f(t')\|_{L^\infty}dt' \leq \frac 1
% {\lam\log2}\quad\hbox{and}
% \quad
\int_0^t 2^{-\Phi_{q,\lam}(t,t')}  (\|S_q\nabla v(t')\|_{L^\infty}+1) dt' \leq
\frac 1 {\lam\log2}\cdotp
$$
Then, we obtain that
\begin{eqnarray*}
M^\s_\lam(f) & \leq  & \|f_0\|_{B^\s_{p,\infty}(H^{-s})}+ \frac C \lam 
M_\lam^{\s+1}(v)
+ C\|\nabla v\|_{\wt L^1_T(C^0)}M^\s_\lam(f) +\frac C \lam
M^\s_\lam(f)\\ 
   & \leq  & \|f_0\|_{B^\s_{p,\infty}(H^{-s})  }+ \frac {C} \lam  M_\lam^{\s+1}(v)
 +\frac {2C} \lam M^\s_\lam(f).
\end{eqnarray*}
This proves the theorem of course if we prove the
estimate (\ref{besovinegapertetransportlemmeeq7}) of the lemma.
First of all, let us decompose the operator~$R_q$. We have
\begin{eqnarray*} 
R_q(v,f) & = & \sum_{\ell=1}^6 R_q^\ell(v,f)\quad {\rm with}\\
R_q^1(v,f) & = & \sum_{j=1}^d\D_q(T_{\partial_j f}v^j),\\
R_q^2(v,f) & = &\sum_{j=1}^d[\D_q, T_{v^j}\partial_j]f,\\
% R_q^3(v,f) & = & \sum_{j=1}^dT_{(v^j-S_qv^j)}\partial_j\D_qf,\\
% R_q^4(v,f) & = & \sum_{j=1}^d-T_{\partial_j\D_qf}S_qv^j,\\
R_q^3(v,f) & = & \sum _{j=1}^d\D_q\partial_jR(v^j,f) + \Delta_{q-1}
 v^j  \partial_j 
\Delta_{q+1} f_q - \Delta_{q-2}  v^j \partial_j 
\Delta_{q-1} f_q  \\
% -R(S_qv^j,\D_q\partial_jf),\\ 
R_q^4(v,f) & = & \sum_{i,j=1}^d \div_g ( c^{ij}_\a  \Delta_q (  T_f  \partial_j v^i)   ) 
   +\div_g (  \Delta_q (  T_f  \nabla_g U       )      ),  \\
R_q^5(v,f) & = & \sum_{i,j=1}^d  \div_g ( c^{ij}_\a  [\Delta_q, T_{\partial_j v^i} ]   f     ) 
  +  \div_g (  [\Delta_q, T_{\nabla_g U } ]   f     )   \\
R_q^6(v,f) & = & \sum_{i,j=1}^d  \div_g \Big(   c^{ij}_\a ( R( \partial_j v^i  ,f) 
 + \Delta_{q-1}   \partial_j  v^i 
\Delta_{q+1} f_q - \Delta_{q-2} \partial_i v^j  
\Delta_{q-1} f_q  )\Big)\\
   &+&  \sum_{i,j=1}^d  \div_g \Big(    R( \nabla_g U   ,f) 
 + \Delta_{q-1}   \nabla_g U 
\Delta_{q+1} f_q - \Delta_{q-2} \nabla_g U   
\Delta_{q-1} f_q   \Big)
%    T_{\partial_j v_i - S_q \partial_j v_i } f_q  )\\
% R_q^9(v,f) & = & - T_{f_q} S_q \partial_j v_i 
% R_q^{10}(v,f) & = & - T_{f_q} S_q \partial_j v_i 
\end{eqnarray*}

 Indeed, 
\begin{eqnarray*}
\D_q(v\cdot\nabla f) & = &
\Delta_q ( \sum_{j=1}^d T_{\partial_j f}v^j + T_{v^j}\partial_j f+ R(v^j,\partial_jf) )
\\
 %  R_q^1(v,f)+
% \sum_{j=1}^d\D_qT_{v^j}\partial_j f+ \D_qR(v^j,\partial_jf)\\
   & = & \sum_{\ell =1}^2R_q^\ell(v,f)+
\sum_{j=1}^d T_{v^j}\partial_j\D_qf+\D_qR(v^j,\partial_jf),\\
%    & = & \sum_{\ell =1}^3R_q^\ell(v,f)+
% \sum_{j=1}^d T_{S_qv^j}\partial_j\D_qf+\D_qR(v^j,\partial_jf).
\end{eqnarray*}
Then, we use that 
\begin{eqnarray*}
\sum_{j=1}^d T_{v^j}\partial_j f_q  &=&   \Sum_{|q-q'|\leq 1 } 
 S_{q'-1}v^j \p_j  \D_{q'} f_q \\
 &=&    S_{q-1}v^j \p_j   f_q  +     \Sum_{|q-q'|\leq 1 } 
 (S_{q'-1}v^j -  S_{q-1}v^j)   \p_j  \D_{q'} f_q \\
 &=&  S_{q-1}v^j \p_j   f_q  +      
 \D_{q-1}v^j    \p_j  \D_{q+1} f_q  -\D_{q-2}v^j    \p_j  \D_{q-1} f_q  
\end{eqnarray*} 
 Hence, 
 \[
\D_q(v\cdot\nabla f) = \sum_{\ell
=1}^3R_q^\ell(v,f)+S_{q-1}v\cdot\nabla
 f_q.
\]
In the same way, we have 
\[
\D_q( \div_g (G( v,f)  f  )  ) = \sum_{\ell
=4}^6R_q^\ell(v,f)+ \div_g (G(S_{q-1} v,S_{q-1}f)  f_q  ). 
\]

Let us estimate the six terms appearing above.  We have

Let us begin with~$R_q^1(v,f)$. By definition of the paraproduct, we have
$$
R_q^1(v,f)=\sum _{j=1}^d\sum_{q'}
\D_q(S_{q'-1}\partial_jf\D_{q'}v^j).
$$
As, if~$|q-q'|> 2$ then the above term is equal to~$0$, we deduce that
$$
\norm {H R_q^1(v(t),f(t))}{L^p(L^2)}  \leq  C\sum_{|q-q'|\leq
2}\|{ H S_{q'-1}\nabla f}\|_{L^\infty(L^2)}\|\D_{q'}v(t)\|_{L^p}.
$$
Using the fact that, if~$|q-q'|\leq 2$, then~$\|{ H S_{q'-1}\nabla
f}\|_{L^\infty(L^2)}\leq C2^q\|H f(t)\|_{L^\infty(L^2)} \leq C 2^q $, 
we infer that
$$
\norm {H R_q^1(v(t),f(t))}{L^p(L^2)} \leq C2^q  \sum_{|q-q'|\leq
2}\|\D_{q'}v(t)\|_{L^p}   \leq C   \sum_{|q-q'|\leq
2}\|\nabla \D_{q'}v(t)\|_{L^p}   .
$$ 
Hence
$$
2^{q\s-\Phi_{q,\lam}(t)}\|H R_q^1(v(t),f(t))\|_{L^p(L^2)}
\leq C  M_\lam^{\s+1}(v) \sum_{|q-q'|\leq 2}
2^{-\lam\int_0^t\|S_q\nabla v(t')\|_{L^\infty}dt'+
\lam\int_0^t\|S_{q'}\nabla v(t')\|_{L^\infty}dt'}.
$$

But, it is obvious that
$$
\int_0^t\|S_{q'}\nabla v(t')\|_{L^\infty}dt'-\int_0^t\|S_q\nabla
v(t')\|_{L^\infty}dt'
\leq \int_0^t\|(S_{q'}-S_q)\nabla v(t')\|_{L^\infty}dt'.
$$
Using the fact that~$|q-q'|\leq 2$,  
we get
\begin{equation}
\label{besovinegapertetransportlemmeeq0}
\int_0^t\|S_{q'}\nabla v(t')\|_{L^\infty}dt'-\int_0^t\|S_q\nabla
v(t')\|_{L^\infty}dt'
\leq C|q-q'|\|\nabla v\|_{\wt L^1_T(C^0)}.
\end{equation}
So it turns out that
\begin{equation}
\label{besovinegapertetransportlemmeeq1}
2^{ q\s-\Phi_{q,\lam}(t)}
\|H R_q^1(v(t),f(t))\|_{L^p(L^2)}
   \leq 2^{C\lam\|\nabla v\|_{\wt L^1_T(C^0)}}
   M_\lam^{\s+1}(v).
\end{equation}

Now let us look at~$R^2_q(v,f)$. By definition of the paraproduct, we have
\begin{eqnarray*}
R_q^2(v,f) = 
 -\sum_{j=1}^d\sum_{q'}[S_{q'-1}v^j\partial_j\D_{q'},\D_q]f\\
     =   -\sum_{j=1}^d\sum_{q'}[S_{q'-1}v^j,\D_q]\partial_j\D_{q'}f.
\end{eqnarray*}
The terms of the above sum are equal to~$0$ except if~$|q-q'|\leq 2$.
Moreover, by
definition of the  operators $\D_q$, we have
$$
[S_{q'-1}v^j,\D_q]\partial_j\D_{q'}f(x)
=2^{qd}\int_{\RR^d}h(2^q(x-y))(S_{q'-1}v^j(x)-S_{q'-1}v^j(y))
\partial_j\D_{q'}f(y)dy.
$$
So we infer that
$$ 
\| H [S_{q'-1}v^j,\D_q]\partial_j\D_{q'}f(x)\|_{L^2(M)}
\leq 2^{-q} \normsup {\nabla
S_{q'-1}v}2^{qd}
\Bigl(\Bigl (2^q|\cdot|\times |h(2^q\cdot)|\Bigr)\star
\|H \partial_j\D_{q'}f\|_{L^2(M)}\Bigr)(x).
$$
Hence, 
$$ 
\| H [S_{q'-1}v^j,\D_q]\partial_j\D_{q'}f(x)\|_{L^p(L^2(M))}
\leq 2^{-q} \normsup {\nabla
S_{q'-1}v} 
\|H \partial_j\D_{q'}f\|_{L^p(L^2(M))}.
$$

Then, we have, using Inequality \refeq{besovinegapertetransportlemmeeq0},
$$ 
\longformule{ 
2^{q\s-\Phi_{q,\lam}(t)}
\|H [S_{q'-1}v^j,\D_q]\partial_j\D_{q'}f\|_{L^p(L^2(M))}
}
{
\leq C M^\s_\lam(f) \sum_{|q-q'|\leq 2}
2^{C\lam \| v\|_{\wt L_T^1(C^1)}}(\|\nabla(S_{q'-1}-S_q)v(t)\|_{L^\infty}
+\|S_qv(t)\|_{L^\infty}).
}
$$
So, we get
\begin{eqnarray} 
2^{q\s-\Phi_{q,\lam}(t)}
\|H R_q^2(v(t),f(t))\|_{L^p(L^2)} \leq C M^\s_\lam(f) 2^{C\lam \| v\|_{\wt
L^1(C^1)}} & &
\nonumber\\
   & &
\label{besovinegapertetransportlemmeeq2}
\!\!\!\!\!\!\!\!\!\!\!\!\!\!\!\!\!\!\!\!\!\!\!\!\!\!\!\!\!\!\!\!\!\!\!\!\!\!
\!\!\!\!\!\!\!\!\!\!\!\!\!\!\!\!\!\!\!\!\!\!\!\!\!\!\!\!\!\!\!\!\!\!\!\!\!\!
\!\!\!\!\!\!\!\!\!\!\!\!\!\!
{}\times\Bigl(\|S_q
\nabla v(t)\|_{L^\infty}+
\sum_{|q-q'|\leq 2} \|\nabla(\D_{q'}v(t)\|_{L^\infty}\Bigr).
\end{eqnarray}

For $R_q^3$, we have 
\begin{eqnarray*}
 \| H {R_q^{3}(v,f)} \|_{L^p(L^2)}   &\leq&  C 
\sum_{\scriptstyle |q'-q''|\leq 1 \atop
\scriptstyle q'\geq q- 2}  2^q
   \| \D_{q'} v   \|_{L^p}    \| H \D_{q''}  f \|_{L^\infty(L^2)}  \\
& \leq &    C 
\sum_{\scriptstyle q'\geq q- 2}  2^{q-q'}
   \| \D_{q'}  \nabla   v  \|_{L^p}    \| H f \|_{L^\infty(L^2)}.  
\end{eqnarray*}
Hence,   
\begin{eqnarray*}
 2^{q\s-\Phi_{q,\lam}(t)}   \| H {R_q^{3}(v,f)} \|_{L^p(L^2)}   &\leq& 
 C 
\sum_{\scriptstyle q'\geq q- 2}  2^{(1+\s)(q-q') - \Phi_{q,\lam}(t) + 
\Phi_{q',\lam} (t) }
         M_\lam^{\s+1}(v)    \| H  f \|_{L^\infty(L^2)}.  
\end{eqnarray*}
Then, we see that the sum converges since
\begin{equation}
|\Phi_{q,\lam}(t) -
\Phi_{q',\lam} (t) | \leq \lam  \|\nabla v\|_{\wt L^1_T(C^0)}   |q-q'| 
\leq (\s -\b) |q-q'|
\end{equation}
and $1+\s -(\s-\b) = 1+\b > 0$. Hence, we get  
\begin{eqnarray*}
 2^{q\s-\Phi_{q,\lam}(t)}   \| H {R_q^{3}(v,f)} \|_{L^p(L^2)}    \leq 
   C          M_\lam^{\s+1}(v)    \| H  f \|_{L^\infty(L^2)}.  
\end{eqnarray*}

The estimate for~$R_q^4(v,f) = R_q^{4,1}(v,f) + R_q^{4,2}(f)  $ 
is the same as the 
estimate for~$R_q^1(v,f)$. Indeed, we have  
\begin{eqnarray*}
\norm {H R_q^{4,1}(v(t),f(t))}{L^p(L^2)}  &\leq&  C\sum_{|q-q'|\leq
2}\|{ \nabla_g H S_{q'-1}  f}\|_{L^\infty(L^2)}\|\D_{q'}\nabla v(t)\|_{L^p} \\
 &\leq&  C\sum_{|q-q'|\leq
2} \|\D_{q'}\nabla v(t)\|_{L^p}
\end{eqnarray*}
where we used that  $ \|{ \nabla_g H S_{q'-1}  f}\|_{L^\infty(L^2)} \leq C  $.
Hence,  we conclude as for~$R_q^1(v,f)$.  Besides, 
\begin{eqnarray*} 
\norm {H R_q^{4,2}(f(t))}{L^p(L^2)}  &\leq&  C\sum_{|q-q'|\leq
2}\|{ \nabla_g H S_{q'-1}  f}\|_{L^\infty(L^2)}\|\D_{q'}\nabla_g U\|_{L^p} \\
 &\leq&  C\sum_{|q-q'|\leq 
2} \|\D_{q'} f(t)\|_{L^p}
\end{eqnarray*}
Hence, we conclude as  for~$R_q^1(v,f)$ and get 
\begin{equation}
\label{besovinegapertetransportlemmeeq5}
2^{ q\s-\Phi_{q,\lam}(t)}
\|H R_q^{5,2}(f(t))\|_{L^p(L^2)}
   \leq 2^{C\lam\|\nabla v\|_{\wt L^1_T(C^0)}}
   M_\lam^{\s}(f).
\end{equation}

We write~$R_q^5(v,f) =  R_q^{5,1}(v,f) + R_q^{5,2}(f)  $.  
The estimate for~$R_q^5(v,f)$ is similar to the one  for~$R_q^2(v,f)$
with the only difference that we have to use the regularity of $\nabla v$. 
 We have
\begin{eqnarray*}
  [\Delta_q, T_{\partial_j v^i} ]   f  = 
%  -\sum_{j=1}^d\sum_{q'}[S_{q'-1}v^j\partial_j\D_{q'},\D_q]f\\
        -\sum_{j=1}^d\sum_{q'}[S_{q'-1} \p_j v^i,\D_q]\partial_j\D_{q'}f.
\end{eqnarray*}
The terms of the above sum are equal to~$0$ except if~$|q-q'|\leq 2$.
Moreover, by
definition of the  operators $\D_q$, we have
$$
[S_{q'-1}\p_j v^i,\D_q] \D_{q'}f(x)
=2^{qd}\int_{\RR^d}h(2^q(x-y))(S_{q'-1}\p_j v^i(x)-S_{q'-1}\p_j v^i(y))
\D_{q'}f(y)dy.
$$
So we infer that
$$   
\| H R_q^{5,1}(v,f) \|_{L^2(M)}
\leq 2^{-q} |\nabla^2
S_{q'-1}v |  2^{qd} 
\Bigl(\Bigl (2^q|\cdot|\times |h(2^q\cdot)|\Bigr)\star
\| \nabla_g  H \D_{q'}f\|_{L^2(M)}\Bigr)(x).
$$
Hence, 
$$ 
\| H  R_q^{5,1}(v,f)  \|_{L^p(L^2(M))}
\leq 2^{-q}   \| \nabla^2
S_{q'-1}v \|_{L^p}  
\| \nabla_g   H \D_{q'}f\|_{L^\infty (L^2(M))}.
$$

Then, we have, using Inequality\refeq{besovinegapertetransportlemmeeq0},
$$ 
\longformule{ 
2^{q\s-\Phi_{q,\lam}(t)}
\|H  R_q^{5,1}(v,f)  \|_{L^p(L^2(M))}
}
{
\leq C  \sum_{|q-q'|\leq 2 \atop q''\leq q'-1}
 2^{(\s-1)(q-q'') - \Phi_{q,\lam}(t) + 
\Phi_{q'',\lam} (t) }   M_\lam^{\s+1}(v)    \| \nabla_g H  \D_{q'} f \|_{L^\infty(L^2)}.
}
$$ 
Hence, 
\begin{eqnarray*}
 2^{q\s-\Phi_{q,\lam}(t)}
\|H  R_q^{5,1}(v,f)  \|_{L^p(L^2(M))} \leq C  \sum_{q''\leq q+1}
2^{-\b(q-q'')} M_\lam^{\s+1}(v)  \| \nabla_g H  f \|_{L^\infty(L^2)} 
\end{eqnarray*}
and the sum is uniformly bounded  
 since $ \s -1 + \lam \|\nabla v\|_{\wt L^1_T(C^0)}   
 \leq - \b   $. Then, we argue in a similar way for 
$ \|H  R_q^{5,2}(f)  \|_{L^p(L^2(M))}  $ and get 
$$ 
\| H  R_q^{5,2}(f)  \|_{L^p(L^2(M))}
\leq 2^{-q}   \| \nabla
S_{q'-1}f \|_{L^p (L^2(M))  }  
\| \nabla_g H  \D_{q'}f\|_{L^\infty (L^2(M))}.
$$
and we conclude as above with $ M_\lam^{\s+1}(v) $ replaced by 
$ M_\lam^{\s}(f) $. 

Finally, the estimate for~$R_q^6(v,f)$ is exactly  the same as the 
 one  for~$R_q^3(v,f)$ since, we also have 
that $\|{ \nabla_g H    f}\|_{L^\infty(L^2)} \leq C$.

\section{Global existence}

\vspace{3mm}
Now, we turn to the proof of our main  theorem.
First, we notice that the local existence in with 
$ v \in L^\infty_{loc} ( [0,T);      W^{1,r} ) \cap L^2_{loc} ( [0,T);  W^{2,r}) $ and 
$ f \in  L^\infty_{loc} ( [0,T);     W^{1,r} (H^{-s})  ) $  can be easily 
deduced from standard arguments.  Moreover, from regularity estimates 
for the heat equation, we have for all  $ 0< T_0 < T$, 
 $v \in  L^\infty_{loc}( (T_0,T); W^{2-\eps,r} )$. 
 
%  We assume
% that we have a solution given by  theorem \ref{existence}
%  on an interval~$[0,T[$. 

We  want to prove that we can extend the solution beyond  the time $T$. 
It is enough to prove that $\nabla v  \in  L^\infty ( (0,T) \times \RR^2)$. 

% Theorem\ref{theoviscoelasticfondsob} 
 
The local existence result tells that, for any~$T_0$
in~$]0,T[$, the solution~$(v,f)$ of~(\ref{NS}) belongs
to the space~$L^\infty_{loc}([T_0,T[;W^{2-\eps,r}\times W^{1,r}(H^{-s}))$
for any $\eps>0$. 
% with $s-\e > {d \over 2}=1$. 
%In the sequel, we take $s$ instead of $s-\e$.
  Sobolev
type embeddings of Corollary \ref{inclusiontypesob} imply that
\[
(v,\tau) \in L^\infty_{loc}\Bigl([T_0,T[;\wt B^{2-\eps-2\left(\frac 1
r-\frac 1 p\right)}_{p,\infty}\times\wt B^{1-2\left(\frac 1
r-\frac 1 p\right)}_{p,\infty}\Bigr).
\]
Choosing $\eps < 1 -2/r$ and~$p=\infty$ in the above assertion implies that~$(v,\tau)
\in L^\infty_{loc}(\wt C^{1+\s} \times \wt C^{\s}(H^{-s})  )$ where 
$\s = 1-\eps - 2/r > 0$.
%   As~$s$ is greater
%than~$1$, the tensor~$\tau$ belongs to~$L^2([T_0,T]; L^2)\cap
%L^1([T_0,T];C^0)$. 
So we can apply 
Theorem \ref{estimation2dprecisee} and  we can  choose~$T_0$ such that,
with the notations of Theorem \ref{besovinegapertetransport}, we
have
\[
\|\nabla v\|_{\wt L^1_{[T_0,T]}(C^0)} \leq   \frac { min(\s-\beta, 1-\s-\b)} {3\lambda}  
 \cdotp
\]
The deteriorating regularity estimate of Theorem \ref{besovinegapertetransport} applied
with~$\s$ and between~$T_0$ and~$T$ tells exactly that~$f$
satisfies
\begin{equation} 
\label{eqprooftheo1200} 
M^{\s}_\lam(f)  \leq  3 \|f\|_{C^{\s} (H^{-s})}+ 
% \biggl(
\frac {3C} \lam
%   + T-T_0\biggr)
M_\lam^{\s+1}(v).
\end{equation}
Now, we have to estimate~$\nabla v$. The two dimensional 
Navier-Stokes equation can
be written as
\[
\partial_t v-\nu \Delta v = P(v\cdot \nabla v) +P D\tau
\] 
where  $P$ denotes the Leray projector on the divergence free vector 
field. Exactly
along the same lines as in the proof of 
Theorem\ref{besovinegapertetransport}, we
have
\[
2^{q(\s+1)-\Phi_{q,\lam}(t)} \|P(v\cdot \nabla v)-P(S_qv\cdot \nabla\Delta _q
v)\|_{L^\infty} \leq C M_\lam^{\s+1}(v) \biggl(\|S_q\nabla
v(t)\|_{L^\infty}+\sum_{q'\geq q}2^{q-q'} \|\nabla\Delta_{q'}
v(t)\|_{L^\infty}\biggr).
\]
Moreover, it is obvious that
\[ 
2^{q(\s-\frac 1 2)-\Phi_{q,\lam}(t)} \|P(S_qv\cdot \nabla\Delta
_q v)\|_{L^\infty} \leq C \|v(t)\|_{H^{\frac 1 2}} M_\lam^{\s+1}(v).
\] 
So it turns out that
\begin{equation}
\label{eqprooftheo1202}
2^{q(\s+1)-\Phi_{q,\lam}(t)} \|\Delta_q P(v\cdot \nabla v)\|_{L^\infty} 
\leq CM_\lam^{\s+1}(v) \biggl(\|S_q\nabla
v(t)\|_{L^\infty}+\sum_{q'\geq q}2^{(q-q')} \|\nabla
v(t)\|_{L^\infty}+2^{\frac {3q} 2}\|v(t)\|_{H^{\frac 1 2}}\biggr).
\end{equation}
Using  well known estimates on the heat equation (see for 
instance\ccite{Chemin99})
and Inequalities\refeq{eqprooftheo1200} and\refeq{eqprooftheo1202} , 
we get that
\[ 
M_\lam^{\s+1}(v) \leq \|v_0\|_{C^{\s+1}}  + \biggl(\frac
C\lam +2^{\frac {3q} 2}F_q(T_0,T)\biggr)M_\lam^{s}(v)  + \frac C\nu
   M_\lam^{\s}(\tau) 
\]
with
\[
F_q(T_0,T) \eqdefa \sup_{t\in[T_0,T]} \int_{T_0}^t e^{c\nu 2^{2q}
(t-t')}\|v(t')\|_{H^{\frac 1 2}}dt'.
\]
H\"older inequality implies immediately that
\[
F_q(T_0,T) \leq \frac C {\nu^{\frac 3 4}}2^{-\frac 
 {3q} 2}\|v\|_{L^4_{T_0,T]}(H^{\frac
1 2})}.
\]
Moreover, it is easy to see that  
\[
M_\lam^{\s}(\tau)  \leq M_\lam^{\s}(f). 
\]
So,  we infer that
\[
M_\lam^{\s+1}(v) \leq \|v_0\|_{C^{\s+1}} 
+ \frac {3C}\nu  \|\tau_0\|_{C^{\s}} + \biggl(\frac
C\lam +  \frac
C{\lam \nu} 
    +\frac C {\nu^{\frac 3 4}} \|v\|_{L^4_{T_0,T]}(H^{\frac
1 2})}\biggr)M_\lam^{\s+1}(v)
\]
Now it is enough to choose~$T_0$ such that the quantity
\[
 \biggl(\frac
C\lam +  \frac
C{\lam \nu} 
    +\frac C {\nu^{\frac 3 4}} \|v\|_{L^4_{T_0,T]}(H^{\frac
1 2})}\biggr)
\]
is small enough. Then as~$\s$ is greater than~$0$, the
solution~$(v,\tau)$ of 
the system (\ref{NS}) is such that~$(\nabla v,\tau)$ belongs
to~$L^\infty([T_0,T]\times\RR^2)$; this concludes the proof
of~Theorem \ref{main}.

\section{Acknowledgments}
The work of P.C. is partially supported by NSF-DMS grant 0504213.
The work of N. M.  is partially supported by NSF-DMS grant  0403983.
% \bibliographystyle{abbrv}
% \bibliography{biblio}

\begin{thebibliography}{10}

\bibitem{BC94}
H.~Bahouri and J.-Y. Chemin.
\newblock \'{E}quations de transport relatives \'a\ des champs de vecteurs
  non-lipschitziens et m\'ecanique des fluides.
\newblock {\em Arch. Rational Mech. Anal.}, 127(2):159--181, 1994.

\bibitem{BCAH87}
R.~B. Bird, C.~Curtiss, R.~Amstrong, and O.~Hassager.
\newblock {\em Dynamics of polymeric liquids, Kinetic Theory Vol. 2,}.
\newblock Wiley, New York, 1987.

\bibitem{Bony81}
J.-M. Bony.
\newblock Calcul symbolique et propagation des singularit\'es pour les
  \'equations aux d\'eriv\'ees partielles non lin\'eaires.
\newblock {\em Ann. Sci. \'Ecole Norm. Sup. (4)}, 14(2):209--246, 1981.

\bibitem{Chemin99}
J.-Y. Chemin.
\newblock Th\'eor\`emes d'unicit\'e pour le syst\`eme de {N}avier-{S}tokes
  tridimensionnel.
\newblock {\em Journal d'Analyse Math\'ematique}, 77(?):27--50, 1999.

\bibitem{CM01}
J.-Y. Chemin and N.~Masmoudi.
\newblock About lifespan of regular solutions of equations related to
  viscoelastic fluids.
\newblock {\em SIAM J. Math. Anal.}, 33(1):84--112 (electronic), 2001.

\bibitem{const} P. Constantin, Nonlinear Fokker-Planck Navier-Stokes
systems, Commun. Math. Sci. {\bf 3} (4), 531-544 (2005).

\bibitem{CFTZ06} P. Constantin, C. Fefferman, E. Titi and A. Zarnescu
\newblock
Regularity for coupled two-dimensional nonlinear Fokker-Planck 
and Navier-Stokes systems  
\newblock {\em Comm. Math. Phys} (to appear) 


\bibitem{doi} M. Doi, S.F. Edwards, {\em The Theory of Polymer Dynamics},
Oxford University Press, 1986.

\bibitem{ELZ04}
W.~E, T.~Li, and P.~Zhang.
\newblock Well-posedness for the dumbbell model of polymeric fluids.
\newblock {\em Comm. Math. Phys.}, 248(2):409--427, 2004.

\bibitem{FGO97}
E.~Fern{\'a}ndez-Cara, F.~Guill{\'e}n, and R.~R. Ortega.
\newblock Some theoretical results for viscoplastic and dilatant fluids with
  variable density.
\newblock {\em Nonlinear Anal.}, 28(6):1079--1100, 1997.

\bibitem{FGO98}
E.~Fern{\'a}ndez-Cara, F.~Guill{\'e}n, and R.~R. Ortega.
\newblock Some theoretical results concerning non-{N}ewtonian fluids of the
  {O}ldroyd kind.
\newblock {\em Ann. Scuola Norm. Sup. Pisa Cl. Sci. (4)}, 26(1):1--29, 1998.

\bibitem{FGOb}
E.~Fern{\'a}ndez-Cara, F.~Guill{\'e}n, and R.~R. Ortega.
\newblock {\em The mathematical analysis of viscoelastic fluids of the
  {O}ldryod kind}.
\newblock 2000.

\bibitem{GS90}
C.~Guillop{\'e} and J.-C. Saut.
\newblock Existence results for the flow of viscoelastic fluids with a
  differential constitutive law.
\newblock {\em Nonlinear Anal.}, 15(9):849--869, 1990.

\bibitem{GS90b}
C.~Guillop{\'e} and J.-C. Saut.
\newblock Global existence and one-dimensional nonlinear stability of shearing
  motions of viscoelastic fluids of {O}ldroyd type.
\newblock {\em RAIRO Mod\'el. Math. Anal. Num\'er.}, 24(3):369--401, 1990.

\bibitem{JLL04}
B.~Jourdain, T.~Leli{\`e}vre, and C.~Le~Bris.
\newblock Existence of solution for a micro-macro model of polymeric fluid: the
  {FENE} model.
\newblock {\em J. Funct. Anal.}, 209(1):162--193, 2004.

\bibitem{LLZ05}
F.-H. Lin, C.~Liu, and P.~Zhang.
\newblock On hydrodynamics of viscoelastic fluids.
\newblock {\em Comm. Pure Appl. Math.}, 58(11):1437--1471, 2005.

\bibitem{LLZ07}
F.-H. Lin, C.~Liu, and P.~Zhang.
\newblock On a Micro-Macro model for polymeric fluids near equilibrium 
\newblock {\em Comm. Pure Appl. Math.}, (to appear)

\bibitem{LM00cam}
P.-L. Lions and N.~Masmoudi.
\newblock Global solutions for some {O}ldroyd models of non-{N}ewtonian flows.
\newblock {\em Chinese Ann. Math. Ser. B}, 21(2):131--146, 2000.

\bibitem{otto} F. Otto, A.E. Tzavaras, Continuity of velocity gradients in suspensions of rod-like molecules, SFB preprint Nr. 141, (2004).

\bibitem{ZZ06arma}
H.~Zhang and P.~Zhang.
\newblock Local existence for the {FENE}-dumbbell model of polymeric fluids.
\newblock {\em Arch. Ration. Mech. Anal.}, 181(2):373--400, 2006.

\end{thebibliography}
%\end{document}

\end{document}